\documentstyle[10pt]{article}
\newlength{\lpr}
\newlength{\wpr}
\input epsf
\def\st{\sqrt{2}}
\def\sth{\sqrt{3}}
\def\qsi{q^{-6}}

\def\qmf{q^{-4}}
\def\qthi{q^{-3}}
\def\qti{q^{-2}}
\def\qmt{q^{-2}}
\def\qi{q^{-1}}
\def\qt{q^{2}}
\def\qth{q^{3}}
\def\qf{q^{4}}
\def\qfi{q^{5}}
\def\qs{q^{6}}
\def\qse{q^{7}}
\def\qe{q^{8}}
\def\qn{q^{9}}
\def\sq{q^{1/2}}
\def\sqi{q^{-1/2}}
\def\tmqs{2_{-q^{6}}}
\def\tmqf{2_{-q^{4}}}
\def\tmqt{2_{-q^{2}}}
\def\tqt{2_{q^{2}}}
\def\tqf{2_{q^{4}}}
\def\tqe{2_{q^{8}}}
\def\thq{3_{q}}
\def\thmq{3_{-q}}
\def\thmqt{3_{-q^{2}}}
\def\thmqf{3_{-q^{4}}}
\def\thqt{3_{q^{2}}}
\def\thqf{3_{q^{4}}}
\def\fqf{4_{q^{4}}}
\def\fq{5_{q}}
\def\fqt{5_{q^2}}
\def\fmq{5_{-q}}
\def\Xp{X_{+}}
\def\Xm{X_{-}}
\def\Xpm{X_{\pm}}
\def\qH{q^{H}}
\def\qHi{q^{-H}}
\def\sqH{q^{H/2}}

\def\sqHi{q^{-H/2}}
\def\eviqr{| i;\, q \rangle}
\def\evjqr{| j;\, q \rangle}

\def\evior{| i;\, 1 \rangle}

\def\eviol{\langle i;\, 1 |}
\def\evjol{\langle j;\, 1 |}
\def\t{{\bf t}}
\def\g{{\bf g}}
\def\l{\lambda}
\def\U{{\cal U}}
\def\F{{\cal F}}
\def\R{{\cal R}}
\def\C{{\cal C}}
\def\P{{\cal P}}
\def\H{{\cal H}}
\def\id{\mathop{\rm id}}
\def\Tr{\mathop{\rm Tr}}
\def\d{{\rm d}}
\def\ot{\otimes}
\newcommand{\Ra}{{\bf R}}
\newcommand{\Rah}{\hat{\Ra}}

\def\tr{\triangleright}

\def\fe{& = &}

\newcommand{\be}{\begin{equation}}
\newcommand{\ee}{\end{equation}}
\newcommand{\bae}{\begin{eqnarray}}
\newcommand{\eae}{\end{eqnarray}}
\def\nn{\nonumber}
\def\ff{\nn \\}

\def\men#1{~(\ref{#1})}
\def\mens#1{~\ref{#1}}

\def\fgr#1#2{\hbox{\raise #1 ex\hbox{\epsfbox{#2}}}}
\def\figr#1#2#3{\hbox{\raise #2 ex\epsfxsize=#1pt
   \hbox{\epsfbox{#3}}}}

\def\bb{\bibitem}
\def\ie{\hbox{\it i.e.\,}}

\begin{document}
\begin{titlepage}
\begin{center}
\vskip 5in
     \hfill   DIAS-STP 97-09\\

\vskip 1.5in

{\large \bf Drinfeld twist for quantum $su(2)$ in the adjoint
representation}

\vskip 1.2in

Chryssomalis Chryssomalakos\footnote[1]{e-mail address:
chryss@stp.dias.ie}\\[.5in]
{\em Dublin Institute for Advanced Studies\\
School of Theoretical Physics\\
10 Burlington Road\\
Dublin 4, Ireland}
\end{center}

\vskip 1.6in

\begin{abstract}
We give a detailed description of the adjoint representation of
Drinfeld's twist element, as well as of its coproduct, for $su_{q}(2)$. 
We also discuss, as applications, the computation of the universal
$R$-matrix in this representation and the problem of symmetrization
of identical-particle states with quantum $su(2)$ symmetry.
\end{abstract}
\end{titlepage}
\renewcommand{\thepage}{\arabic{page}}
\setcounter{page}{1}
\section{Introduction}

Drinfeld's work on quasitriangular quasi-Hopf algebras appeared 
(in english) in 1990. Its impact on elucidating the conceptual
foundations of the deformation of semi-simple Lie algebras cannot
be overestimated. By relaxing coassociativity, one perceives,
through Drinfeld's work,  
both the classical and the deformed algebras to belong to the same
``space'', with a universal twist element $\F$ effecting the 
``rotation'' 
of one into the other. However, despite an existence of proof, this
twist element has proved particularly elusive, the only known
explicit result being for the case of the $q$-deformed Heisenberg
algebra $\H_{q}(1)$\cite{BGST}. For the particularly important case
of $sl_{q}(2)$ (from which $\H_{q}(1)$ can be obtained by
contraction) reference~\cite{CGZ} provides the ``semi-universal''
expression $(\rho_{1/2} \ot \id) \F$ while~\cite{DNS} contains $\F$
up to the second order in $h \equiv \ln q$. 
This lack of complete information on $\F$ has been however less 
unsatisfactory a situation than one might expect since, quite
often, all that is needed in applications is a 
matrix representation of the twist. 
 A number of results exist in this direction
- reference~\cite{Engel} provides such representations in the
tensor square of the fundamental for $su(N)$, $so(N)$ and 
$sp(N)$ while~\cite{Cu,Za} contain  expressions for 
$(\rho_{1/2} \ot \rho_{j})\F$, where $\rho_{j}$ denotes the
$(2j+1)$-dimensional representation of $sl_{q}(2)$.

There are two results of particular interest not included in the
above list. On the one hand, the representation of $\F$ in the
tensor square of the adjoint of $sl_{q}(2)$ is not known, while on
the other, no information exists in the literature on
representations of $(\Delta \ot \id) \F$, a matrix which solves the
problem of a symmetrization procedure compatible with quantum group
actions. Both computations involve rather tedious algebra, due to
the size of the matrices involved, yet explicit results are
desirable in applications. We undertake therefore in this paper an
explicit computation of the above quantities, pointing out along
the way arguments that reduce the complexity of the task. 

The layout of the paper is as follows. In sections 2 and 3 we
briefly review the key ingredients that enter in the problem,
namely the adjoint representation of $su_{q}(2)$, using the
Jimbo-Drinfeld basis, and the Drinfeld-Kohno theorem respectively
(we assume $q$ to be real throughout the paper, so that $\F$ is
represented by an orthogonal matrix - this is the reason for
referring to $su_{q}(2)$ above rather than to $sl_{q}(2)$). Section
4 contains the explicit computations mentioned above along with a
detailed description of the structure of the resulting matrices.
Applications appear in section 5, where the adjoint representation 
of the universal $R$-matrix is computed and a quantum permutation
operator is constructed which commutes with the quantum symmetry 
generators.
\section{Quantum $su(2)$ in the adjoint representation}

Quantum $su(2)$ (or, more precisely, the quantum deformation of the
universal enveloping algebra of $su(2)$, denoted by $\U_{q}(su(2))
\equiv \U_{q}$) is the Hopf algebra generated by 
$q^{\pm H/2}$, $\Xpm$ obeying the commutation relations
\bae
q^{H/2} \Xpm \fe q^{\pm 1} \Xpm q^{H/2} \ff
\Xp \Xm - \Xm \Xp \fe \l^{-1} (\qH - \qHi) \, , \label{su2alg}
\eae 
where $\l = q - \qi$ and $q$ is a real number (the latter
requirement originates in the $*$-structure of the algebra which we
do not discuss here).
 These tend, in the classical limit (\ie as
$q$ approaches 1), to the familiar $su(2)$ commutation relations.
The coproduct is given by
\bae
\Delta(q^{\pm H/2}) \fe q^{\pm H/2} \ot q^{\pm H/2} \ff
\Delta(\Xpm) \fe \Xpm \ot q^{H/2} + q^{- H/2} \ot \Xpm \label{suq2cop}
\eae 
and is an algebra homomorphism with the (standard)
multiplication $ (x \ot y) (z \ot w) = xz \ot yw $ in the tensor 
square of the algebra.
The counit $\epsilon$ and antipode $S$ are given by
\be
\epsilon(q^{\pm H/2})=1, \qquad \epsilon(\Xpm)=0, \qquad 
S(q^{\pm H/2})= q^{\mp H/2}, \qquad S(\Xpm) = -q^{\pm 1} \Xpm \, .
\ee
We examine now the adjoint representation of the above algebra. 
\subsection{The factorized basis}
We first recall the definition of the adjoint representation in the
quantum case. Classically, this representation is defined via the
adjoint action of the generators among themselves, which closes
linearly in the generators. The representation is then extended to
the entire universal enveloping algebra as a homomorphism.
 In the quantum case, the appropriate definition for the adjoint
action of $x$ on $y$ ($x$, $y \in \U_{q}$) is given by
\be
x \tr y = x_{(1)} y S(x_{(2)})
\ee
which reduces to a commutator for the classical coproduct
$\Delta(x) = x \ot 1 + 1 \ot x$. One easily checks that this
adjoint action does not close linearly in the set $\{q^{H/2}, \Xpm \}$.
It is easily seen though that it does close linearly for the set of
generators $x_{1}, x_{+}, x_{-}, x_{2}$, given by
\bae
x_{1} \fe \l^{-1}(1-\qHi -\qi \l^{2} \Xp \Xm) \ff
x_{+} \fe \sqi \Xp q^{H/2} \ff
x_{-} \fe \sqi q^{H/2} \Xm \ff
x_{2} \fe \l^{-1}(1-\qH) \, , \label{xopmt} 
\eae  
thus giving rise to the folowing representation
\be
\rho_{f} \left( \begin{array}{c|c}
x_{1} & x_{+} \\ \hline
x_{-} & x_{2}
\end{array} \right) =
\left( \begin{array}{cccc|cccc}

-q^{-2} \l & 0 & 0 & \l & 0 & 0 & \qi & 0 \\
0 & \qi +q^{-3} -q & 0 & 0 & -\qi & 0 & 0 & q \\
0 & 0 & -q & 0 & 0 & 0 & 0 & 0 \\
q^{-2} \l & 0 & 0 & -\l & 0 & 0 & -\qi & 0 \\ \hline
0 & -\qi & 0 & 0 & 0 & 0 & 0 & 0 \\
0 & 0 & 0 & 0 & 0 & -q & 0 & 0 \\
q^{-3} & 0 & 0 & -\qi & 0 & 0 & \qi & 0 \\
0 & \qi & 0 & 0 & 0 & 0 & 0 & 0
\end{array} \right) \, . 
\ee

Inverting now\men{xopmt} we find
{\small
\be
\rho_{f}(q^{H/2}) \! = \! \left( \! \! \begin{array}{cccc}
1 &  0 &  0 &  0 \\
0 &  q &  0 &  0 \\
0 &  0 &  \qi &  0 \\
0 &  0 &  0 &  1
\end{array} \! \! \right) \! ,  
\rho_{f}(\Xp) \! = \! \left(\! \! \! \begin{array}{cccc}
0 & 0 & \sq & 0 \\
-\sqi & 0 & 0 & q^{3/2} \\
0 & 0 & 0 & 0 \\
0 & 0 & -\sq & 0
\end{array} \! \! \! \right) \! ,  
\rho_{f}(\Xm) \! = \! \left( \! \! \begin{array}{cccc}
0 & -\sqi & 0 & 0 \\
0 & 0 & 0 & 0 \\
q^{-3/2} & 0 & 0 & -\sq \\
0 & \sqi & 0 & 0
\end{array} \! \! \right) \!  .
\ee
}
\subsection{The reduced basis}
\label{Trb}
To simplify subsequent calculations we chose to work in a reduced
basis, in which the representation (denoted $\rho_{r}$) consists 
of 1 and 3-dimensional 
blocks. The transition is effected by conjugation with the matrix 
$A$, \ie ($x \in \U_{q}$)
\be
\rho_{f}(x) = A \rho_{r}(x) A^{-1}
\ee
with $A$ being given by
\be
A = \left( \begin{array}{cccc}
q^{3} & 0 & \gamma^{-1} & 0 \\
0 & 1 & 0 & 0 \\
0 & 0 & 0 & \qi \\
q & 0 & -\gamma^{-1} & 0
\end{array} \right) \, .
\ee
This amounts to switching to the generators 
\be
c = q^{3} x_{1} + q x_{2} \, , \qquad x_{0} = \gamma^{-1}(x_{1} -
x_{2})
\ee
in favor of $x_{1}$, $x_{2}$ as well as rescaling $x_{-}$ by $q$ (the
order chosen is \{$c$, $x_{+}$, $x_{0}$, $x_{-}$\}). The
generator $c$ is central and provides therefore an 1-dimensional
representation for the  algebra, given by the counit (notice that
this coincides with the classical limit). Omitting the
1-dimensional block, the reduced representation is given explicitly
by
\be
\rho_{r}(\Xp) \! = \! \left( \! \begin{array}{ccc}
0 & -\sqi \gamma & 0 \\
0 & 0 & \sqi \gamma \\
0 & 0 & 0
\end{array} \!  \right) ,  \, \, \,
\rho_{r}(q^{H/2}) \! = \! \left( \! \begin{array}{ccc}
q & 0 & 0 \\
0 & 1 & 0 \\
0 & 0 & \qi
\end{array} \! \right) , \, \, \,
\rho_{r}(\Xm) \! = \! \left( \! \begin{array}{ccc}
0 & 0 & 0 \\
-\sqi \gamma & 0 & 0 \\
0 & \sqi \gamma & 0 
\end{array} \! \right) , 
\label{rhor}
\ee
where $\gamma \equiv \sqrt{1+\qt}$. The rescaling of $x_{-}$ 
mentioned above guarantees that
$\rho_{r}(\Xp)^{T} = \rho_{r}(\Xm)$, as in the classical case. 
It is also worth pointing out that, under the substitution $q
\mapsto \qi$, $\rho_{r}(\Xpm)$ go into themselves while $\rho_{r}(\qH)
\mapsto \rho_{r}(\qH)^{-1}$. 
We introduce now, for later
use, the quantum casimir $\C_{q}$, given by
\be
\C_{q}= \Xm \Xp + \l^{-2}
(q^{H+1} + q^{-(H+1)} - 2) \, .
\ee
Its value in $\rho_{r}$ is ${[3/2]}^{2}$ where $[x] \equiv
\frac{q^{x}-q^{-x}}{q-\qi}$. Notice that $\lim_{q \rightarrow
1} \C_{q} = \C_{c}+1/4$, where $\C_{c}$ is the familiar undeformed
casimir of $su(2)$.
\section{The Drinfeld-Kohno theorem}

We give a brief discussion of the essential points only, as they
pertain to the particular case examined (in particular, we do not
mention coassociativity since it is not essential to the problem at
hand) - the
reader is referred to~\cite{Drin} for a detailed account of
these matters in a general context.

 Starting from the (undeformed) 
universal enveloping algebra $\U(su(2)) \equiv \U$ one can form 
the algebra
$\U([[h]])$ by extending the field of coefficients to the ring of
formal power series in $h$.
 With the identification
$h \equiv \ln(q)$, $\U([[h]])$ can be shown to be isomorphic
to $\U_{q}$. In other words, one can find $h$-dependent
invertible functions of the classical $su(2)$ generators, which obey the
deformed relations\men{su2alg}. Applying the classical coproduct to
these functions one does not however obtain the coproduct
of\men{suq2cop} (with $H$, $\Xpm$ in these equations considered as 
functions of the classical generators). Instead, as the
Drinfeld-Kohno theorem states, the following relation holds
\be
\Delta(x) = \F \Delta_{c}(x) \F^{-1} \label{FDFi}
\ee
for $x \in \U_{q}$, where $\F \equiv \F^{(1)} \ot \F^{(2)}
\in \U_{q}^{\ot 2}$ is called the {\em twist}. In the above relation 
$\Delta(\cdot)$
stands for the coproduct that appears in\men{suq2cop} while
$\Delta_{c}(\cdot)$ is the classical cocommutative coproduct given
by $\Delta_{c}(y)=y \ot 1 + 1 \ot y$ where $y$ is any classical Lie
algebra generator.  

As an illustration of the content of the theorem, consider the
particular case where $x$ in\men{FDFi} stands for a classical
generator (then, in the rhs, $\Delta_{c}(x)$ does not involve $q$).
To compute the lhs, one could express $x$ in terms of the quantum
generators $\qH, \, \Xpm$ (via the isomorphism mentioned in the
theorem),  then take the coproduct
using\men{suq2cop} and finally switch back to the classical
generators using the inverse of the above isomorphism. Notice that 
the $q$-dependence of the lhs comes,
in this case, entirely from the argument of the coproduct 
($x$ being a $q$-dependent
function of the quantum generators) - the coproduct\men{suq2cop},
written in terms of $K \equiv \qH, \, \Xpm$ does not involve $q$
explicitly. On the rhs, the $q$-dependence comes entirely through $\F$.
 One way to look at\men{FDFi} is to consider it as defining a
second ($q$-dependent) coproduct for the classical algebra $\U$ via
conjugation by $\F$.

There is more to the Drinfeld-Kohno theorem though. $\U$ is a
cocommutative coassociative quasitriangular Hopf algebra and as such its
universal $R$-matrix is trivial (\ie equal to $1 \ot 1$). It can
however {\em alternatively} be considered as a quasitriangular 
quasi-Hopf algebra with universal $R$-matrix given by $\R_{c}=q^{\t}$
where 
\be
\t = \Delta_{c}(\C_{c})-\C_{c} \ot 1 - 1 \ot \C_{c} \, .
\label{tdef}
\ee
Adopting this point of view one sees the universal
$R$-matrix of $\U_{q}$ to be the image, under the twist, of
its classical counterpart according to
\be
\R = \F' \R_{c} \F^{-1} 
\label{RFpRcFi}
\ee
where $\F' \equiv \P (\F)$ and $\P$ is the permutation
operator in $\U_{q}^{\ot 2}$. 
\section{Adjoint representation of the twist}
\label{Arott}

We compute here explicitly the adjoint representation of $\F$,
following (a variation of) the method presented in~\cite{Engel}. 
\subsection{Preliminaries}
\subsubsection*{The method}
Consider the matrix representation $C_{q} \equiv 
\rho_{r}^{\ot 2}(\Delta(\C_{q}))$
of the coproduct of the quantum casimir. Its eigenvectors $\eviqr$,
$i=1 \ldots 9$, form (for each value of $q$) a basis in the nine 
dimensional space on which the
above matrix acts ($i$ can be thought of as a composite label,
$i=\{ \alpha_{q}, \, h_{q} \}$, where $\alpha_{q}$ is an 
eigenvalue of the casimir in $\rho_{r}^{\ot 2}$ and $h_{q}$ 
differentiates between
eigenvectors belonging to the same casimir eigenvalue - we will use
the eigenvalues of $H$ for that purpose). Taking representations 
in both spaces of\men{FDFi} (with $\C$ replacing $x$) one finds
\be
C_{q} = F_{q}  \rho_{r}^{\ot 2}
 (\Delta_{c}(\C_{q})) F_{q}^{-1}
\ee
where $F_{q} \equiv \rho_{r}^{\ot 2}(\F)$. It follows that
\be
\eviqr = F_{q} \evior ,  \qquad i=1 \ldots 9 \, .
\label{evqFevo}
\ee
Putting $M_{ij} \equiv \evjol i; \, 1 \rangle$ we find for $F_{q}$
\be
F_{q} = M^{-1}_{ij} \evjqr \eviol  \, .
\label{FMev}
\ee 
Notice that we have allowed for the possibility of the eigenvectors
being nonorthogonal. Trying to produce an orthonormal set by taking
suitable linear combinations might (and will, typically) involve
mixing eigenvectors of different (casimir) eigenvalues.
Nevertheless, in the particular representation chosen 
($\rho_{r}^{\ot 2}$), the casimir will turn out to be given by a
symmetric matrix and, consequently, the matrix $M$ above will be
diagonal, resulting in significant simplification of the
computations. 
\subsubsection*{Orthogonality}
One of the advantages in switching from $\rho_{f}$ to $\rho_{r}$ is
the fact that, as mentioned earlier, in this latter representation,
the transpose of $\Xp$ is equal to $\Xm$ (this requirement dictated
the rescaling of $x_{-}$). A glance at the form of
the coproduct of $\C_{q}$ (equation\men{Cqcop} below) reveals that
the above property (along with $H$ being diagonal) results in
$C_{q}$ being a symmetric matrix. Its
eigenvectors therefore can be taken to form, for each value of $q$, 
an orthonormal basis. Equation\men{evqFevo} then shows that $F_{q}$ is
an orthogonal matrix which rotates this basis from its classical
(at $q=1$) to its quantum  position (at general $q$).
\subsubsection*{Invariant subspaces}
The coproduct of $H$ being undeformed, $F_{q}$ commutes
with its representation and is therefore of a
block-diagonal form, with each block acting in a subspace of fixed
$H$-eigenvalue. We expect therefore two 1-dimensional blocks
(corresponding to the directions $(++)$ and $(--)$ with $H$-eigenvalues
4 and -4 respectively), two 2-dimensional blocks (corresponding to
the planes $(+0,0+)$ and $(0-,-0)$ with $H$-eigenvalues 2 and -2
respectively) and one 3-dimensional block (acting in the subspace
$(+-,00,-+)$ of $H$-eigenvalue 0). Similar remarks hold for
$C_{q}$. The 1-dimensional blocks of $F_{q}$ are equal to 1 due to
orthogonality while those of $C_{q}$ are 
obviously equal to each other (being the values of the casimir 
in the same
multiplet). A glance at eqn.\men{Cqcop} below shows that the
2-dimensional blocks of $C_{q}$ are also equal to each other (since
$H$ has opposite eigenvalues in the corresponding planes). It
then follows, from \men{FMev}, that the above property
holds for $F_{q}$ as well. 
We will present any matrix $W$, with the structure described above,
 in the form $\{ \mu_{W}, \, D_{W}, \, T_{W} \}$ where $\mu_{W}$, 
$D_{W}$, $T_{W}$ are 1, 2 and 3-dimensional matrices respectively.
\subsection{Explicit computation of $F_{q}$}
We start by computing $C_{q}$. We first find
\bae
\Delta(\C_{q}) \fe \Xm \Xp \ot \qH + \Xm \sqHi \ot \sqH \Xp \ff
  & &                 +\sqHi \Xp \ot \Xm \sqH + \qHi \ot \Xm \Xp \ff
  & &              + \l^{-2} (q \cdot \qH \ot \qH +\qi \qHi \ot \qHi
                - 2 \cdot 1 \ot 1) 
\label{Cqcop}
\eae
and substituting from\men{rhor} we obtain
\be
C_{q}=\{ {[5/2]}^{2}, \, 
\left( \begin{array}{cc}
{[1/2]}^{2}+\qi {[2]}^{2} & [2] \\ 
{[2]}^{} & {[1/2]}^{2} + q{[2]}^{2}
\end{array} \right) \, ,
\left( \begin{array}{ccc}
{[1/2]}^{2} + q^{-2}[2] & -\qi [2] & 0 \\
-\qi [2] & {[1/2]}^{2} + 2[2] & -q[2] \\
0 & -q[2] & {[1/2]}^{2} + q^{2} [2]
\end{array} \right)
 \} \, . 
\ee
The corresponding 2 and 3-dimensional eigenvectors are
\be
\begin{array}[b]{cc}
(1, \, -q^{-2})   & {[3/2]}^{2} \\
(1, \, q^{2})  & {[5/2]}^{2} \\
(q^{2}, \, q, \, 1)  & {[1/2]}^{2} \\
(1, \, -\l, \, -1) &  {[3/2]}^{2} \\
(1, \, -q -q^{3}, \, q^{4})  & {[5/2]}^{2} 
\end{array} \, ,
\label{Cqev}
\ee
where, in the second column, we give the corresponding eigenvalue
($\C_{q}$ has eigenvalues of the form ${[j+1/2]}^{2}, \, j=0, 1, 2,
\ldots $). It is now easy to construct $F_{q}$ using\men{FMev}. We
find
\be
\mu_{F_{q}}= 1, \qquad D_{F_{q}}= 
\frac{1}{\sqrt{2}a} \left( \begin{array}{cc}
q^{2}+1 & -(q^{2}-1) \\
q^{2}-1 & q^{2}+1
\end{array} \right)  
\ee
as well as
\be
T_{F_{q}} = \frac{1}{\sqrt{6}ab} \left( \begin{array}{ccc}
\st \qt a +\sth q b +1 & \st \qt a -2 & \st \qt a-\sth qb+1 \\
\st qa -\sth (\qt-1)b -q(1+\qt) & \st qa+2q(1+\qt) & \st qa +\sth
(\qt -1)b -q(1+\qt) \\
\st a -\sth qb +q^{4} & \st a-2q^{4} & \st a+\sth qb +q^{4}
\end{array} \right) \,, 
\label{T}
\ee
where $a \equiv \sqrt{1+q^{4}}, \, b \equiv \sqrt{1+\qt +q^{4}}$.
Recalling the discussion in section\mens{Trb}, regarding the
properties of $\rho_{r}$ under the substitution $q \mapsto \qi$, 
we see that, with $\Delta(\C_{q})$ as in\men{Cqcop}, $C_{q}$
satisfies ${(C_{\qi})}_{12} = {(C_{q})}_{21}$. It follows then
from\men{FMev} that the same relation holds for $F_{q}$.
\subsubsection*{The 2-dimensional subspace}
$D_{F_{q}}$ is an $SO(2)$ matrix, even in $q$, corresponding 
to an angle of rotation $\theta^{(2)}_{q}$ given by
\be
\tan(\theta^{(2)}_{q})=\frac{\qt -1}{\qt +1} \, ,
\ee
counterclockwise in the $(+0,\, 0+)$ or $(0-, \, -0)$ plane. As
$q$ ranges in the interval $(0, \, \infty )$, $\theta^{(2)}_{q}$
goes from $-\frac{\pi}{4}$ to $\frac{\pi}{4}$. Notice that although
at the endpoints of the above interval in $q$, neither the quantum
$su(2)$ algebra nor the adjoint representation\men{rhor} make sense,
$D_{F_{q}}$ retains nevertheless a well defined value.  Under $q
\mapsto \qi$, the two axes get interchanged and the angle of
rotation changes sign.

As pointed out in~\cite{Engel}, a more natural reference point in
describing the rotation effected by $F_{q}$, and one which
simplifies the results, is the {\em crystal
limit} $q \rightarrow 0$ (notice that the eigenvectors of $C_{q}$
in\men{Cqev}, when normalized, tend to the basis vectors in this
limit). The matrix that connects the general
value of $q$ to this limit is $F_{q0} \equiv F_{q} F_{0}^{-1}$. We
find
\be
D_{F_{q0}} = \frac{1}{a} \left( \begin{array}{cc}
1 & -\qt \\
\qt & 1
\end{array} \right) \, ,
\ee
with corresponding angle of rotation $\theta^{(2)}_{q0}$ given by
$\tan(\theta^{(2)}_{q0})=\qt$.

\subsubsection*{ The 3-dimensional subspace}
$T_{F_{q}}$ is an $SO(3)$ matrix - we denote by $\hat{n}_{q}$ the
unit vector along the axis of the rotation represented by $T_{F_{q}}$
and by $\theta^{(3)}_{q}$ the angle of rotation (clockwise) around 
$\hat{n}_{q}$
($\theta^{(3)}_{q}$ ranges in the interval $[0, \, \pi)$). Working
in a basis where one of the basis vectors is along $\hat{n}_{q}$,
and using the invariance of the trace of $T_{F_{q}}$ under rotations, 
we find
\be
1+2\cos(\theta^{(3)}_{q})= \Tr (T_{F_{q}}) \, , 
\label{thetaTrT}
\ee
where, using\men{T},
\be
\Tr(T_{F_{q}}) = \frac{1}{\sqrt{6}ab} (\st(\qt +q+1)a +2\sth qb
+q^{4} +2q^{3} +2q+1) \, .
\label{TrT}
\ee
Let $\alpha_{i}, \, i=1,2,3$ denote the angles that $\hat{n}_{q}$
makes with the axes - a little geometry shows that
\be
\cos^{2}(\alpha_{i}) =
\frac{{T_{F_{q}}}_{ii}-\cos(\theta^{(3)}_{q})} 
     {1-\cos(\theta^{(3)}_{q})} \, .
\ee
Notice that ${T_{F_{q}}}_{ii}$ (no summation over $i$) is the angle 
that the $i$-th basis
vector makes with its image under the rotation. The correct signs
for $\cos(\alpha_{i})$ are obtained by examining the derivatives $\d
F_{q}/\d q |_{q=0}$ and $\d^{2} F_{q}/\d q^{2}|_{q=0}$. Under $q \mapsto
\qi$, the axes $(+-)$ and $(-+)$ get interchanged. As $q$ ranges
over the reals, $T_{F_{q}}$ traces out a curve in the $SO(3)$
manifold (a ball of radius $\pi$ with antipodal points on the
surface identified).
We have plotted  this curve (for the range $-1\leq q \leq 1$) in 
figure\mens{TFq}.
\begin{figure}
\epsfbox{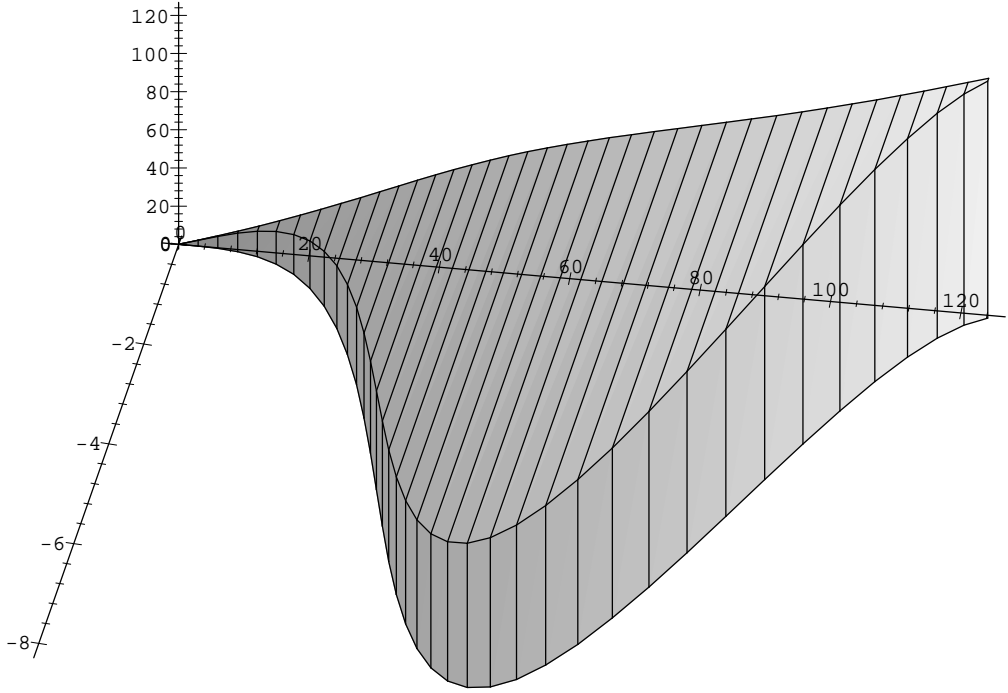}
\caption{The curve $T_{F_{q}}$ in the $SO(3)$ manifold 
($-1 \leq q \leq 1$)}
\label{TFq}
\end{figure}
The left - right axis is in the $(+-)$-direction while the
vertical is along $(-+)$ (the axes are marked in degrees). The points 
of the curve are projected
down to the $(+-,00)$-plane and back to the
$(+-,-+)$-plane (at equal intervals in $q$) for clarity
of presentation. 
To find the value $q_{s}$ at which the curve
meets the surface of the ball, we substitute\men{TrT}
in\men{thetaTrT}, set $\theta^{(3)}_{q}$ equal to
$\pi$ and solve for $q_{s}$ - the result is $q_{s}=-1$. The
corresponding rotation matrix has $(1,-1,1)$ along the second
diagonal as its only nonzero entries. 
The curve starts (for
$q=-1$) at the right part of the figure (upper edge of the
ribbon), touching at that point the $(+-,-+)$-plane on the diagonal
(notice the anisotropic scaling of the axes) as well as the surface
of the $SO(3)$ ball. 
At $q=1$, $T_{F_{q}}$ is the unit matrix, which corresponds to the
origin of the ball. 
Notice that the tangent to the
curve at the origin is along the diagonal in the $(+-,-+)$-plane
(the same holds at $q=-1$). 
The property ${(F_{\qi})}_{12} = {(F_{q})}_{21}$ reflects itself in the
following symmetry of the curve: 
the points corresponding to $q$ and $\qi$ 
 are mapped into each other by central
reflection and subsequent interchange of the $(+-)$, $(-+)$ axes
(this allows one to picture the missing segments $-\infty <q<-1$
and $1 < q < \infty$,
which are omitted from figure\mens{TFq} for clarity).
Notice that the above two missing segments ``touch'' as $q \rightarrow
\pm \infty$ since, in this limit, 
$T_{F_{q}}$ tends to the common value
\be
T_{F_{\pm \infty}} = \frac{1}{\sqrt{6}} \left( \begin{array}{ccc}
\st & \st & \st \\
-\sth & 0 & \sth \\
1 & -2 & 1
\end{array} \right) \, .
\ee

As in the 2-dimensional case, the results are simplified when the
crystal limit is used as reference point. We find
\be
T_{F_{q0}} = \frac{1}{ab} \left( \begin{array}{ccc}
1 & qb & \qt a \\
-q(1+\qt) & (1-\qt)b & qa \\
q^{4} & -qb & a
\end{array} \right) \, .
\ee
The corresponding $SO(3)$ curve (for $q$ in the interval $[-3,
3]$) is plotted in figure\mens{TFq0}. In this case, the lower
left - upper right axis is along the $(+-)$ direction while the
vertical is again along $(-+)$.  
\begin{figure}
\epsfbox{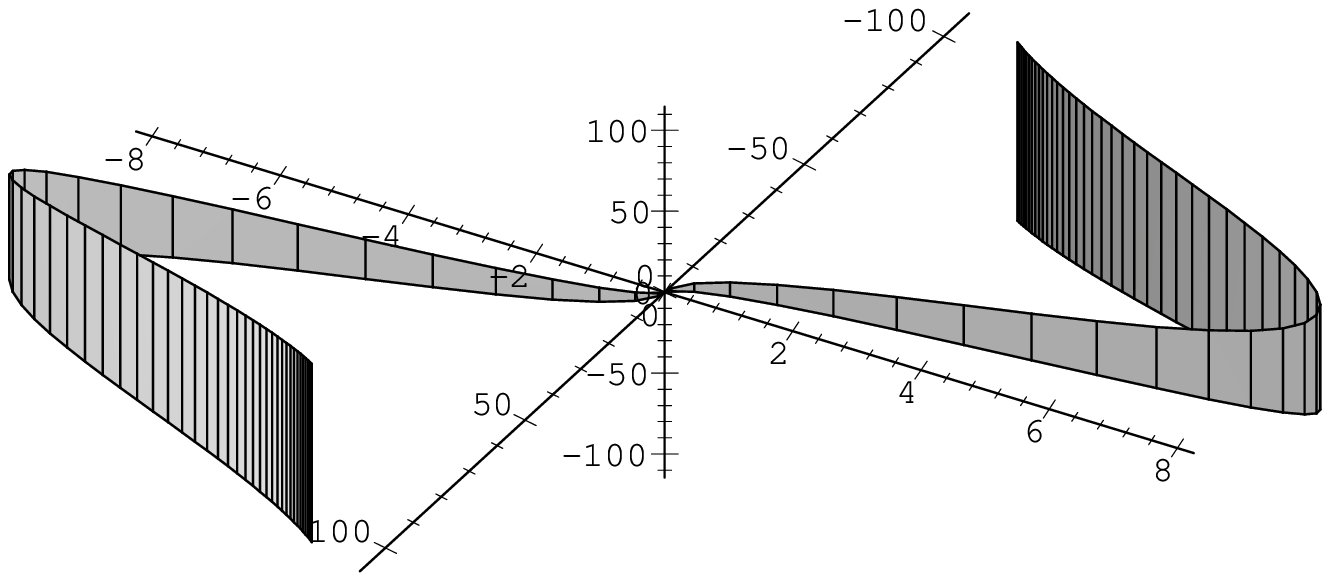}
\caption{The curve $T_{F_{q0}}$ in the $SO(3)$ manifold 
($-3 \leq q \leq 3$)}
\label{TFq0}
\end{figure}
The origin corresponds to $q=0$ while the point of
the curve on the surface of the ball (not shown) to $q \rightarrow 
\pm \infty$.

A remark on the omission of $c$ from
$\rho_{r}$ is in order at this point. It is sometimes desirable 
to know $F_{q}$ in the entire
16-dimensional space spanned by the bilinears in $c, \, x_{+}, \,
x_{0}, \, x_{-}$. This, for example, could originate in the fact
that although the quantum adjoint action closes linearly in the set
$x_{+}, \, x_{0}, \, x_{-}$, the {\em algebra} of these three generators
involves $c$. Such cases are easily handled noting that $F_{q}$ is
simply a unit matrix in the 7-dimensional space we have omitted in the
above computation.
\subsection{The representation of the coproduct of $\F$}
The method for computing $F_{q}$ used above is easily adapted to
the computation of ${F_{q}}_{123} \equiv \rho_{r}^{\ot 3}((\Delta \ot
\id ) \F)$. The matrix $C_{123}=\rho_{r}^{\ot 3}(\Delta^{(2)}  \C)$,
where $\Delta^{(2)}(\C) \equiv (\Delta \ot \id) \circ
\Delta(\C)$, 
is again
symmetric and its eigenvectors are rotated by the orthogonal 
${F_{q}}_{123}$ from their classical to their deformed position. This
allows the computation of ${F_{q}}_{123}$ from a formula analogous
to\men{FMev}. Both $C_{123}$ and ${F_{q}}_{123}$ are of block-diagonal
form, with each block acting in a subspace of fixed $H$-eigenvalue.
Based on this observation, one easily sees that both matrices have 
1, 3 and
6-dimensional blocks (two of each) and one 7-dimensional one. Blocks 
of equal
dimension are equal to each other, up to (easily identified) 
permutations of rows and
columns. We have carried out, with the help of MAPLE, the explicit
computation of these matrices, with respect to the crystal limit
(\ie we have computed the matrix ${F_{q0}}_{123} \equiv
{F_{q}}_{123}{F_{0}}_{123}^{-1})$. We present  the (rather cumbersome)
results in the appendix.
\section{Applications}
\subsection{The universal $R$-matrix in the adjoint representation}
The block diagonal form of the matrices involved simplifies
significantly the computation of $\Ra \equiv \rho^{\ot 2}_{r}(\R)$.
Applying $\rho^{\ot 2}_{r}$ in\men{RFpRcFi} gives
\be
\Ra = F_{q}' \, \rho^{\ot 2}_{r}(q^{\t}) \, F_{q}^{-1} \, .
\label{RaFRcFi}
\ee
For the representation of the classical $R$-matrix we compute
\be
\rho_{r}^{\ot 2}(\t) = \{ 2, \, \left( \begin{array}{cc}
0 & 2\\
2 & 0
\end{array} \right) \, , 
\left( \begin{array}{ccc}
-2 & -2 & 0 \\
-2 & 0 & -2 \\
0 & -2 & -2
\end{array} \right) \}
\ee
which gives, upon exponentiation
\be
\mu_{R_{c}}= \qt, \qquad D_{R_{c}}= \frac{1}{2\qt}
\left( \begin{array}{cc}
\qf +1 & \qf -1 \\
\qf -1 & \qf +1
\end{array} \right)
\ee
as well as
\be
T_{R_{c}} = \frac{1}{6\qf} \left( \begin{array}{ccc}
\qs +3\qt +2 & -2(\qs-1) & \qs -3\qt +2 \\
-2(\qs -1) & 2(2\qs +1) & -2(\qs -1) \\
\qs -3\qt +2 & -2(\qs -1) & \qs +3\qt +2
\end{array} \right) \, ,
\ee
where $R_{c} \equiv \rho_{r}^{\ot2}(\R_{c})$. We may now
use\men{RaFRcFi} to find
\be
\mu_{\Ra} = \qt, \qquad D_{\Ra} = \left( \begin{array}{cc}
1 & \qt -\qti \\
0 & 1
\end{array} \right) , \qquad
T_{\Ra} = \left( \begin{array}{ccc}
\qti &  -q +\qthi &  (\qt-\qti)(1-\qti)^{2} \\
0 & 1 & -q+\qthi  \\
0 & 0 & \qti
\end{array} \right)  .
\ee
$\Ra$ is seen to be upper triangular while $\Rah \equiv P \Ra$
is symmetric ($P$ stands here for the permutation matrix,
$P_{ij,kl}= \delta_{il} \delta_{kj}$ - $\Rah$ is obtained from
$\Ra$ by interchange of the rows of $D_{\Ra}$ and interchange of
the first and third row of $T_{\Ra}$). One could also compute $\Ra$
starting from the known expression for $\R$, apply $\rho_{f}^{\ot 2}$ to
it and then conjugate with the tensor square of $A$ - the two
results agree. 

The characteristic equation for $\Rah$ is easily computed now.
Indeed, one finds that
\be
(D_{\Rah} -\qt I) (D_{\Rah} +\qti I) = 0,
\qquad (T_{\Rah} -\qt I) (T_{\Rah} +\qti I) (T_{\Rah} -\qmf I) =0
\, .
\ee
Notice that $\mu_{\Rah}$ and $D_{\Rah}$ both satisfy the
characteristic equation for $T_{\Rah}$ - this implies that $\Rah$
does too, \ie $\Rah$'s minimal polynomial is cubic. The
characteristic polynomial for $\Rah$ is obtained from its minimal
one by raising its factors, in the order they appear above, to the 
powers 5, 3 and 1. Had we included
$c$ in $\rho_{r}$, obtaining in this way a 16-dimensional matrix
for $\Rah$, the minimal polynomial would turn out to be quintic,
the extra two factors being $(\Rah -I)(\Rah +I)$, coming from the
$c$-sector in which $\Rah$ is a permutation matrix. The above two
factors would be raised to the power 3 and 4 respectively in the
characteristic polynomial of $\Rah$.
\subsection{Quantum symmetrization}
We consider here
the problem of quantum symmetrization of identical particle states. 
We supply, for reasons of self-containment, a brief exposition of 
the relevant theory, along the lines of~\cite{FS}, which contains a
thorough discussion.

As a first step towards the implementation of quantum group
symmetries in quantum mechanics and quantum field theory, the
question of the compatibility of Bose/Fermi statistics with quantum
group actions ought to be addressed. Consider, for concreteness, a
two (identical) particle system in quantum mechanics. The Hilbert
space of the system is split into symmetric and antisymmetric
subspaces, a division which is respected by the action of classical
Lie algebra generators. This latter fact is traced to the symmetry
of the classical coproduct under exchange of the two spaces, so
that the classical permutation matrix $P$ 
commutes with the representation of the
coproduct of the generators
\be
[P,\rho^{\ot 2}(\Delta_{c}(x))]=0, \qquad x \in U(\g),
\ee
where $\g$ is the Lie algebra under consideration and $\rho$ the
representation through which it acts on the one-particle states. 
In the quantum case,
the lack of cocommutativity implies that the quantum group action
mixes the two eigenspaces of $P$. The problem is rectified with the
introduction of a quantum permutation operator $P^{F}$, given by 
conjugation
by $F$ of the classical one, \ie \, $P^{F}=FPF^{-1}$. One easily
shows that 
\be
[P^{F}, \, \rho^{\ot 2}(\Delta(x))]=0 \, \qquad x \in U_{q}(\g) \, ,
\ee
which guarrantees that quantum symmetric and antisymmetric states
(defined as eigenstates of $P^{F}$)
remain such after being acted upon by elements of $U_{q}(\g)$
(in the above equation, $F$ is computed in $\rho$). 

Specifying $\g$, in the discussion above, to be 
$su(2)$ and working in the representation $\rho_{r}$, we find 
\be
S^{F} = \{ 1, \qquad \frac{1}{1+\qf} \left( \begin{array}{cc}
1 & \qt \\
\qt & \qf
\end{array} \right) , \qquad 
\frac{1}{1+\qf} \left( \begin{array}{ccc}
\qf -\qt +1 & \qth -q & \qt \\
\qth -q & \qt & -\qth +q \\
\qt & -\qth +q & \qf -\qt +1
\end{array} \right) \}
\ee
while
\be
A^{F} = \{ 0, \qquad \frac{1}{1+\qf} \left( \begin{array}{cc}
\qf & -\qt \\
-\qt & 1
\end{array} \right) , \qquad
\frac{1}{1+\qf} \left( \begin{array}{ccc}
\qt & -\qth +q & -\qt \\
-\qth +q & (\qt -1)^{2} & \qth -q \\
-\qt & \qth -q & \qt
\end{array} \right) \}
\ee
where $S^{F} \equiv F \frac{1}{2}(1+P) F^{-1}$ and $A^{F} \equiv F 
\frac{1}{2}(1-P)F^{-1}$ are the quantum symmetrizer and
antisymmetrizer respectively. 
In the $c$-sector, $S^{F}$ and $A^{F}$ are equal to their classical
counterparts. Clearly, ${S^{F}}^{2} =
S^{F}$ (similarly for $A^{F}$) and  $S^{F} A^{F}=0$.

For systems with more than two identical particles one proceeds
along similar lines. The quantum permutation matrix for adjacent spaces
$i$, $i+1$ is given by $P^{F}_{i,i+1}={F_{q}}_{12\ldots n} P_{i,i+1}
{F_{q}}^{-1}_{12 \ldots n}$, where $P_{i,i+1}$ is the corresponding
 classical permutation matrix and ${F_{q}}_{12\ldots n}$ is the
representation of $({\Delta}^{(n-1)} \ot \id ) \F$ (different choices
for ${F_{q}}_{12\ldots n}$ are possible). 
For example, for $n=3$, ${F_{q}}_{123}$ can be obtained from the
crystal limit 
result of section\mens{Arott} as ${F_{q}}_{123} = {F_{q0}}_{123}
{F_{10}}^{-1}_{123}$.

\section*{Appendix}
We give here the explicit form 
for the 3, 6 and 7-dimensional blocks of ${F_{q0}}_{123}$, acting in 
the subspaces
$(++0,+0+,0++)$, $(++-,+00,+-+,0+0,00+,-++)$ and
$(+0-,+-0,0+-,000,0-+,-+0,-0+)$ respectively. 
{\small
\be
{F^{(3)}_{q0}}_{123}  =  \frac{1}{mk} \left( \begin{array}{ccc}
m & \qf k & -\qs \\
-\qf m & k & -\qt \\
0 & \qt k & k^{2}
\end{array} \right) \, , \nn
\ee
\vskip .5cm
\be
{F^{(6)}_{q0}}_{123}  =  \left( \begin{array}{cccccc}
y & -\qth w & \qf z & 0 & \qe r & \qt st^{2} \\
\qth \thqt x^{2}y & \tmqs w & -\qi \tqt z & -\qt x & -\qth \tqt r &
q \tqf^{-1} s \\
0 & \qth w & \qmf z & 0 & r & s \\
\qfi (1+\thqt) x^{2}y & \qt w & -q \tqt z & x & (q-\qse +\qn)r &
-\qi \tqt st^{2} \\
-\qf \thqt x^{2} y & -\qf w & -q^{-3} \tqt z & \qf x & (\qi -q
+\qse) r & (q\tqf)^{-1} s \\
(-\qf +\qs +q^{10}) x^{2}y & 0 & z & q \tqt x & - \tqe r & \qmt
st^{2}
\end{array} \right) \, , \nn
\ee
\vskip .5cm
while ${F^{(7)}_{q0}}_{123}$ is given by
\be
\left( \begin{array}{ccccccc}
\thqt^{-1/2} & \qt \thmqf u & -\qt t & (-1+\qf + \qs) d & e & (-\qf
-\qs + q^{10}) f & g \\
0 & (1-\qt +\qs + \qe) u & \qt t & \tmqf \thqt d & \qmf e & \qf
\tqt^{2} f & 0 \\
0 & (\qf -\qs -\qe) u & t & -\thqt d & \qt e & (\qt -\qs -\qe) f &
g \\
q \thqt^{-1/2} & (-\qt +\qf +\qs) \tmqf u & -q \tmqt t & \qi (2-\fqt)
d & -\qthi \tqt e & \qth (2+\qt -\qf -\qs +\qe) f &
-\qthi \tmqf g \\
0 & (-\qt +\qf +\qs) u & -\qf t & -\qt \thqt d & \qsi e & (\qf -
\qe -q^{10}) f & -\qti g \\
0 & -2\qs u & -\qt t & 0 & e & \tqt \tqe f & \qti \tmqf g \\
\qt \thqt^{-1/2} & (-\qt -\qf +\qs) u & \qt t & d & \qmf e & -\qt
\thmqf f & (2-\thqt) \qmf g
\end{array} \right) \! \! . \nn
\ee
}
We have used the notation
$n_{x} \equiv 1+x+x^{2} + \ldots +x^{n-1}$ and
\be
\begin{array}{rclcrcl}
d \fe q(8_{\qt} \thqt)^{-1/2} & \qquad \qquad &
e \fe \qs (\tqf \thmqt \fmq \fq)^{-1/2} \ff
f \fe (\tqt \tqf \tqe \fmq \fq)^{-1/2} & \qquad \qquad &
g \fe \qf ((1-\qt +\qs +\qe) \tqf)^{-1/2} \ff
k \fe \sqrt{\tqe} & \qquad \qquad &
m \fe \sqrt{\thqf} \ff
r \fe q(\tqf \tqe \thmq \thq \thmqt)^{-1/2} & \qquad \qquad &
s \fe \qt \sqrt{\frac{\thqt}{\fqt}} \ff
t \fe (\fqt +\qf)^{-1/2} & \qquad \qquad &
u \fe ((1-\qt +\qs +\qe)\tqf \thq \thmq \thmqt)^{-1/2} \ff
w \fe \fqf^{-1/2} & \qquad \qquad &
x \fe (\fqt +2\qf)^{-1/2} \ff
y \fe \tqf^{2} \thq \thmq x^2 & \qquad \qquad &
z \fe \qf (\thq \thmq \thmqt \fmq \fq)^{-1/2} \, . \nn
\end{array}
\ee
\end{document}